\magnification=\magstep1
\hfuzz 11pt
\baselineskip 17.5 pt
\hoffset 1cm
\hsize 340pt
\vsize 500pt

\null
\vskip 1cm
\centerline {\bf IRREDUCIBLE REPRESENTATIONS } 
\vskip .3cm
\centerline {\bf OF THE QUANTUM ANALOGUE OF SU(2)} 
\vskip 2.62cm 
\centerline {P. Roche and D. Arnaudon}
\vskip .5cm
\centerline {Centre de Physique Th{\'e}orique}
\centerline {Ecole Polytechnique F91128 Palaiseau Cedex}
\vskip 2.62cm
\par
{\it Abstract:} \par
\vskip .4cm
We give the complete set of irreducible representations of ${\cal U}(SU(2))_q$ 
when $q$ 
is a $m-th$ root of unity. In particular we show that their dimensions are less 
or equal to $m$. Some of them are not highest weight representations.
\vfil
\line {A854.0988 \hfil September 1988}
\eject

\null
\par
{\bf I.} Introduction \par
\vskip .5cm
Quantum groups take their origin in the context of completely integrable 
systems, where they have shown to be a powerful tool to construct $R$ matrices 
[Drinfeld, Jimbo]. ${\cal U}(SU(2))_q$ 
is the simplest case of this algebraic structure, it is the trigonometric limit 
of the Sklyanin algebra [Sklyanin]. 
It is the symmetry group of the $XXZ$ Heisenberg chain [Pasquier-Saleur] 
and it can be 
used to construct completely integrable systems as well as conformal field 
theory [Pasquier]. Among the values of $q$, the roots of unity correspond to 
discrete series in conformal field theories.\par
In this letter we show that: (let $q$ be a $m-th$ root of unity)
\vskip .2cm
\item {+} {for odd $m$, ${\cal U}(SU(2))_q$ has $N-$dimensional irreducible 
representations if $1\leq N \leq m$.}
\item {+} {for even $m$, ${\cal U}(SU(2))_q$ has $N-$dimensional irreducible 
representations if $1\leq N \leq {m\over 2}$ or $N=m$.}
\par In both cases, $m$ dimensional representations are not always highest 
weight representations.
\vskip .4cm
\par We will give the complete set of irreducible representations of 
${\cal U}(SU(2))_q$. In particular, we prove that the set of finite dimensional 
irreducible unitary representations of ${\cal U}(SU(2))_q$ given by Sklyanin is 
complete. 
\vskip .4cm
\par
In section II we recall the results for $q$ not a root of unity and give some 
preliminaries. In section III, we classify the representations of dimension 
strictly less than $m$ and prove that they are highest weight representations. 
In section IV we prove that there is no irreducible representation of dimension 
strictly greater than $m$. In section V we study the $m-
$dimensional representations and prove that some of them are not 
of highest weight.
\vskip 1.5cm
\eject
{\bf II.} Preliminaries \par
\vskip .5cm
Let $q\in {\bf C}^{\star}-\{1,-1\}$. ${\cal U}(SU(2))_q$ 
is the Hopf algebra over ${\bf C}$ generated by 
$q^{H/2}$, $q^{-H/2}$, $J^+$, $J^-$ satisfying the relations:
$$q^{H/2}q^{-H/2}=q^{-H/2}q^{H/2}=1$$
$$q^{H/2}J^{\pm}q^{-H/2}=q^{\pm}J^{\pm}$$
$$[J^+,J^-]={q^{H}-q^{-H} \over q-q^{-1}}=(H)_q$$
We will not use the coalgebra structure in the following.\par
\vskip .3cm
The link between Sklyanin algebra (in the trigonometric case) and 
${\cal U}(SU(2))_q$ is as follows:
$$S_0={q-q^{-1} \over 2i} . {q^{H/2}+q^{-H/2} \over q^{1/2}+q^{-1/2}}$$
$$S_1=-{q-q^{-1} \over 2i} . {q^{H/2}-q^{-H/2} \over q^{1/2}-q^{-1/2}}$$
$$S_+={q-q^{-1} \over 2i} J^-$$
$$S_-={q-q^{-1} \over 2i} J^+$$
(in Sklyanin's notation, $e^{2i\eta}=q$.)
\vskip .2cm \par
Recall that if $q$ is not a root of unity, the finite dimensional 
representations of ${\cal U}(SU(2))_q$ 
are completely classified by $(N,\omega) \in {\bf 
N}^{\star}\times \{1,-1,i,-i\}$ [Rosso, Lusztig]. 
They are highest weight representations. 
One can find a basis $(w_p)_{0\leq p\leq N-1}$ such that:
$$\cases   {q^{H/2} w_p &$= \omega q^{{N-1 \over 2}-p}w_p $ \cr
            J^+ w_p     &$= \omega [(p)_q(N-p)_q]^{1 \over 2} w_{p-1} 
                           \hskip .7cm  1\leq p\leq N-1$\cr
            J^+ w_0     &$= 0$ \cr
            J^- w_p     &$= \omega [(p+1)_q(N-p-1)_q]^{1 \over 2} w_{p+1}
                          \hskip .7cm  0\leq p\leq N-2$\cr
            J^- w_{N-1} &$= 0 $\cr}\leqno (1)$$
where $(\alpha)_q={{q^{\alpha}-q^{-\alpha}}\over {q-q^{-1}}}$.
We use the same determination of the square root for $J^+$ and $J^-$. 
In this basis, $J^+=(J^-)^t$.
\vskip .8cm
\par We suppose in the following 
that $m$ is the smallest integer such that $q^m=1$.
\vskip .3cm \par
Let us begin with the following remark: $(J^+)^m$, $(J^-)^m$ and $(q^{H/2})^m$ 
are in the center of the algebra. (This is a direct consequence of the 
commutation relations.)
\vskip .3cm \par
{\bf Lemma 1:}{\it Let $M$ be a finite dimensional simple module over {\bf C}. 
Then 
$q^{H/2}$ is diagonalisable and $M=\oplus _{p\in {\bf Z}_m} M_{\lambda q^p} $, 
where $M_{\mu}$ is the eigenspace of $q^{H/2}$ associated to the eigenvalue 
$\mu$.}
\par \vskip.2cm
{\bf Proof:} 
\item{a)} {${\cal U}(SU(2))_q$ is spanned as a vector space by $(J^-)^r 
(J^+)^s (q^{H/2})^t$ where $r,s,t \in {\bf N}^2 \times {\bf Z}$ }
\item{b)} {Let $v$ be an eigenvector of $q^{H/2}$ associated to $\lambda$. Then
$M={\cal U}(SU(2))_q.v$ since $M$ is simple. Because $(J^-)^r (J^+)^s .v$ is an 
eigenvector of $q^{H/2}$ the lemma follows.}
\vskip 1.2cm
\par
{\bf III.} Classification of simple modules of dimension $N< m$.
\vskip .5cm
{\bf Theorem:}\par     
{\it 
\item {a)} {If $m$ is odd and if $1\leq N< m$ the module is defined by the 
relations (1)}
\item {b)} {If $m$ is even and if $1\leq N< {m \over 2}$ 
the module is defined by the relations (1) }
\item {c)} {If $m$ is even, representations of dimension ${m \over 2}$ are 
labelled by $\lambda \in {\bf C}$.}
\item {d)} {If $m$ is even, there is no irreducible representation of 
dimension $N$ for ${m \over 2}<N<m$.}
All these representations are highest weight representations.}
\par\vskip.3cm
{\bf Proof:} When $N< m$ lemma 1 implies that the representation is a 
highest weight representation (same proof as in [Rosso]). As a result there 
exists a basis $(v_p)_{0\leq p\leq N-1}$ such that:
$$\cases   {q^{H/2} v_p &$= q^{{\mu \over 2}-p}v_p  $ \cr
            J^+ v_p     &$= [(p)_q (\mu+1-p)_q] v_{p-1} 
               \hskip .5cm  1\leq p\leq N-1$ \cr
            J^+ v_0     &$= 0 $ \cr
            J^- v_p     &$= v_{p+1} 
               \hskip .5cm  0\leq p\leq N-2$ \cr
            J^- v_{N-1} &$= 0 $ \cr} $$
We have ($\forall p \in [1,N]$ $(p)_q\not= 0$) $\Longleftrightarrow$ ($m$ even 
and $N<{m\over 2}$, or $m$ odd)). Consequently in these cases we must have
$[(N)_q(\mu+1-N)_q]=0$, which implies that $q^{\mu/2}=\omega q^{(N-
1)/2}$. a) and b) are thus proved. 
When $m$ is even and $N\geq {m\over 2}$ these 
relations define a simple module only for $N={m\over 2}$, and $\mu$ is then a 
free parameter. Two modules with different values of $q^{\mu/2}$ are not 
isomorphic.
After a suitable renormalization of the basis,
$$\cases   {q^{H/2} w_p &$= q^{{\mu \over 2}-p}w_p  $ \cr
            J^+ w_p     &$= [(p)_q(\mu+1-p)_q]^{1 \over 2} w_{p-1} 
               \hskip .5cm  1\leq p\leq N-1$ \cr
            J^+ w_0     &$= 0 $ \cr
            J^- w_p     &$= [(p+1)_q(\mu-p)_q]^{1 \over 2} w_{p+1} 
               \hskip .5cm  0\leq p\leq N-2$ \cr
            J^- w_{N-1} &$= 0 $ \cr} \leqno (2)$$
using the same determination of the square root for $J^+$ and $J^-$. 
In this basis, $J^+=(J^-)^t$. 
\vskip 1.2cm
\par
{\bf IV.} Classification of simple modules of dimension $N> m$.
\vskip .5cm
\par Let $N\geq m$
\vskip .3cm
\par
{\bf First case:} $J^+$ and $J^-$ are injections (in that case the 
representation is not a highest weight representation).
\vskip .3cm\par
{\bf Proposition 1:} $N=m$
\vskip .2cm\par
Proof: since $J^-$ is an injection, all the $M_{\lambda q^p}$ have the same 
dimension $n$, and hence $N=mn$. Let $v_0^{(i)},i=1,\dots,n$ be a basis of 
$M_{\lambda}$. Define
$$v_p^{(i)}=(J^-)^p v_0^{(i)} 
\hskip 1.5cm 1\leq i\leq n \hskip 1.5cm 0\leq p\leq 
m-1$$
We must have
$$\cases   {q^{H/2} v_p^{(i)} &$= q^{{\mu \over 2}-p}v_p^{(i)}  $ \cr
            J^-     v_p^{(i)} &$= v_{p+1}^{(i)} 
               \hskip .5cm  0\leq p\leq N-2$ \cr
            J^-v_{m-1}^{(i)}  &$= \sum_j a^i_jv_0^{(j)} $ \cr
            J^+v_0^{(i)}      &$= \sum_j b^i_jv_{m-1}^{(j)} $ \cr
                }  \leqno (3)$$
These relations define two matrices $A$ and $B$. The relation 
$$[J^+,(J^-)^p]=(J^-)^{p-1}(p)_q(H-p+1)_q$$
allows us to calculate 
$$J^+v_p^{(i)}=(p)_q(\mu-p+1)_q v_{p-1}^{(i)} + \sum_{j,k} b^i_j a^j_k v_{p-
1}^{(k)}$$
An easy calculation shows that these relations define a module if and only if 
$[A,B]=0$. Since $^tA$ and $^tB$ commute they have a common eigenvector 
$(x^1,\dots,x^n)=X$ such that $^tAX=ax$ and $^tBX=bx$. Define $V_p=\sum_i 
x^i v_p^{(i)}$ for 
$0\leq p\leq m-1$. $\{V_p\}_{0\leq p\leq m-1}$ is a submodule. 
Henceforth $N=m$.
\vskip .3cm
\par
{\bf Second case:} $J^+$ or $J^-$ is not injective. For example we can suppose 
it exits $v\not=0$, $v\in M_{\lambda}$ such that $J^+v=0$.
\vskip .3cm\par
{\bf Proposition 2:} $N=m$
\vskip .2cm\par
Proof: $N=nm+r$ with $0\leq r <m$. Suppose $r\not= 0$. Then $w,$ $J^-w,\dots$,
$(J^-)^{N-1}w$ is a basis of the module and
$$(J^-)^N w= \sum_{0\leq k\leq n-1}\alpha_k (J^-)^{km+r} w$$
Then $(J^-)^r$ is not an injection. Let $w_0\in M_{\lambda}-\{0\}$ 
such that $(J^-)^rw_0=0$. Since $J^+(J^-)^{km}w=0$ for $k=0,\dots,n$, 
$Vect\{w_0,J^-w_0,\dots,(J^-)^{r-1}w_0\}$ is a non-zero 
submodule of $M$. As a result 
$r=0$, i.e. $N=nm$. This is a particular case of the {\bf first case} with 
$B=0$. The same conclusion leads to $N=m$.
\vskip .2cm\par
We have thus proved that there are no simple module of dimension greater than 
$m$.
\vskip 1.2cm
\par
{\bf V.} Irreducible representations of dimension $m$.
\vskip .5cm
In this case one can construct a basis $(v_p)_{0\leq p\leq m-1}$ of the module 
such that:
$$\cases   {q^{H/2} v_p &$= q^{{\mu \over 2}-p}v_p  $ \cr
            J^- v_p     &$= [(p+1)_q(\mu-p)_q+\alpha \beta ]^{1 \over 2} 
                             v_{p+1} 
               \hskip .5cm  0\leq p\leq m-2$ \cr
            J^- v_{m-1} &$= \alpha v_0 $ \cr
            J^+ v_p     &$= [(p)_q(\mu+1-p)_q+\alpha \beta ]^{1 \over 2} 
                             v_{p-1} 
               \hskip .5cm  1\leq p\leq m-1$ \cr
            J^+ v_0     &$= \beta v_{m-1} $ \cr
         } \leqno (4)$$
where $\mu$ is chosen such that $[(p)_q(\mu+1-p)_q+\alpha 
\beta ]\not=0$ $\forall p\in 
\{1,\dots,m-1\}$. This module is an highest weight module in the case $\alpha 
\beta =0$. 
It is always irreducible unless $m$ is even and $\alpha \beta =0$. \par
Note that $J^+=(J^-)^t$ when $\alpha =\beta $, and $J^+=(J^-)^+$ when 
$\alpha =\beta^{\star}$ and $\mu$ is real. In this last case, our 
representations with three real parameters
correspond to those obtained by Sklyanin. 
\par
Let us denote $M(\mu,\alpha ,\beta )$ this module. 
Two non highest weight representations
$M(\mu,\alpha ,\beta )$ and $M(\mu',\alpha ',\beta ')$ 
are isomorphic if, and only if
$$\mu'=\mu+2r \hskip 1cm r\in {Z}$$
$${\alpha ' \over \beta '} = {\alpha  \over \beta }$$
$$\alpha \beta -\alpha '\beta '= (2r)_q (\mu+2r+1)_q$$
\vskip 1.2cm
\par
{\bf VI.} Conclusion
\vskip .5cm
The finite dimensional irreducible representations of ${\cal U}(SU(2))_q$ 
are now 
classified for every value of $q$.
\par\vskip .2cm
The quantum analog of Hermann-Weyl theorem on complete reducibility of semi-
simple Lie algebras proved by M. Rosso in the case $q$ not a root of 
unity is not valid when $q^m=1$.
\par For example, let $M$ be the module defined by $(v_p)_{0\leq p\leq m} $
and the following relations:
$$q^{H/2} v_p = q^{{m \over 2}-p}v_p$$
$$J^-v_p=v_{p+1}$$
$$J^-v_m=J^+v_0=0$$
$$J^+v_p=(p)_q(1-p)_q v_{p-1}$$
${\bf C}v_m$ is a submodule of $M$, but $M$ is not completely reducible.
\par\vskip .3cm
Our work is just a step toward the understanding of representations of 
${\cal U}(SU(2))_q$. 
A number of questions are left unanswered:
\par - Classification of non semi-simple representations of ${\cal U}(SU(2))_q$.
\par - Generalization of our work to other quantum Lie algebras.
\par - A better understanding of the connection between representations of 
quantum groups and conformal field theories.
\par - It would also be interesting to understand the role (if any) of the 
whole set of representations of ${\cal U}(SU(2))_q$ in the construction of 
solutions of the Yang-Baxter equation.
\vskip 1cm
We thank the referee for having brought to our attention the important work of 
Sklyanin. 
\vskip 1.5cm
\par
{\bf REFERENCES:}
\vskip 1cm
\item {$\star$ }
{V. G. Drinfeld. Quantum groups. Proc. I. C. M. Berkeley. (1986)}
\vskip .25cm
\item {$\star$ }
{M. Jimbo. Quantum R matrix for the generalized Toda system. C. M. P. 
102, 537-547 (1986)}
\vskip .25cm
\item {$\star$ }
{E. K. Sklyanin. Some algebraic structures connected with the Yang-Baxter 
equation. Representations of quantum algebras. Functional Analysis and its 
Applications 17, (1983) p273.}
\vskip .25cm
\item {$\star$ }{V. Pasquier and H. Saleur. Seminar Les Houches 1988}
\vskip .25cm
\item {$\star$ }
{V. Pasquier. Continuum limit of lattice models built on quantum 
groups. Saclay Preprint SPhT/87-125. }
\vskip .25cm
\item {$\star$ }
{V. Pasquier. Etiology of IRF models. CMP 118, (1988) p355. }
\vskip .25cm
\item {$\star$ }
{M. Rosso .Finite dimensional representations of the quantum analogue of 
the enveloping algebra of complex simple Lie algebras. C. M. P. 117,581 (1988)}
\vskip .25cm
\item {$\star$ }
{G. Lusztig. Advances in Mathematics 70, (1988) p237.}

\bye